\documentclass[a4paper,12pt]{article}
\usepackage{amssymb}
 \usepackage{amsthm}
\usepackage{latexsym}
\usepackage{amsmath}
\usepackage{amssymb}
\usepackage{color}

\newcommand{\udots}{\mathinner{\mskip1mu\raise1pt\vbox{\kern7pt\hbox{.}}
\mskip2mu\raise4pt\hbox{.}\mskip2mu\raise7pt\hbox{.}\mskip1mu}}
\usepackage{amsmath}

\newtheorem{theorem}{Theorem}[section]
\newtheorem{corollary}{Corollary}[section]
\newtheorem{lemma}{Lemma}[section]
\newtheorem{remark}{Remark}[section]

\begin{document}
\title{\bf \Large {\bf   Stochastic representations and probabilistic characteristics of multivariate  skew-elliptical distributions}}
{{\author{\normalsize{Chuancun Yin}
\thanks{Corresponding author.}
\\{\normalsize\it    (School of Statistics and Data Science, Qufu Normal University}\\
\noindent{\normalsize\it Shandong 273165, China}\\
email:  ccyin@qfnu.edu.cn)\\
\normalsize{ Narayanaswamy Balakrishnan}
\\
{\normalsize\it (Department of Mathematics and Statistics, McMaster University,}\\
\noindent{\normalsize\it Hamilton, Ontario, Canada}\\
  email: bala@mcmaster.ca)\\
}
\maketitle

\noindent{\large {\bf Abstract}}  The family of multivariate skew-normal distributions has many interesting properties. It is shown
here that these hold for a general class of skew-elliptical distributions.  For this class, several stochastic representations are established and   then their probabilistic properties, such as  characteristic function, moments,  quadratic forms  as well as transformation properties, are  investigated.   }



\baselineskip =16pt

\numberwithin{equation}{section}
\section{Introduction}\label{intro}

In many applied fields such as finance and economics,   data sets  often do not follow the normal
distribution, and sometimes error structures in a regression-type model do not  satisfy symmetry property. In many instances, there is a presence
of high skewness. To model departures from normality, researchers have been active  in constructing   flexible
parametric classes of multivariate distributions  that  possess   skewness and kurtosis and thus     differ from  the
normal distribution. The skew-normal distribution is   one such important model that generalizes the normal distribution to allow for skewness, and has become a widely-used parametric model in multivariate data analysis, due to its tractability and many appealing properties. The pioneering works of
Azzalini (1985, 1986) laid  the foundations for the  skew-normal distribution. The multivariate version of the skew-normal distribution was introduced by Azzalini and
Dalla Valle (1996) and Azzalini and Capitanio (1999). In the past two decades, the class of multivariate skew normal distributions has been studied by many authors, and   the associated developments have been   well-documented by Genton (2004), Lee and McLachlan (2022)  and Azzalini (2022).

In recent years,  there has been considerable interest in the subject of   skew-elliptical distributions which  facilitate the simultaneous modelling of  asymmetry and tail weight behavior. While an initial formulation of   skew-elliptical distributions was mentioned  by   Azzalini and  Capitanio (1999), Branco and Dey (2001)  considered its construction leading to the multivariate  skew-normal distribution  based on  conditioning mechanism. Azzalini and Capitanio (2003)  then studied    skew-elliptical distributions and   showed  that, like the  skew-normal  family, they can be obtained by an additive construction approach. Some
basic properties and the joint distribution of several quadratic forms are  obtained by Fang (2003, 2005).  Kim and Genton (2011) derived the characteristic function of the multivariate  skew-normal  distribution and of its scale mixtures.

As in the case of  multivariate skew-normal distributions, there are other equivalent ways of defining  skew-elliptical distributions through different stochastic representations. These representations may be used for finding moments and also for proposing simulational algorithm.  Branco and Dey (2001)  and Azzalini and Capitanio (2003) obtained  multivariate  skew-elliptical distribution   by conditioning on one random variable being positive, while   Azzalini and Capitanio (2003) derived skew elliptical densities through transformation method and   then demonstrated that the two stochastic representations are equivalent.  A more detailed description of these and other skewed models may be found in the works
 of  Genton (2004),  Azzalini and  Capitanio (2014), Lee and McLachlan (2022), and Azzalini (2022).

This work adds  to the existing literature in several directions. Firstly, we present some  more stochastic representations for skew-elliptical random vectors.
  These representations are  direct generalizations of those of elliptical distributions and have many applications including the derivation of  moments of skew-elliptical distributions via the moments  of skew-normal distributions.
 Secondly, we derive the characteristic function of multivariate skew-elliptical distribution which is useful in examining
 stochastic orderings (e.g.,  Amiri and Balakrishnan (2022) and  Amiri et al. (2020)).  We specifically provide  a negative answer to a  conjecture made by  Genton and   Loperfido (2005).  Thirdly, we derive the  moments and  also of quadratic forms  for multivariate skew-elliptical distributions which are closely  related to the  multivariate skewness and kurtosis measures (e.g.,  Mardia (1970),  Kollo (2008), and Balakrishnan and Scarpa (2012)).
Many of these have  been established   for some specific subclasses such as skew-normal distributions and scale-mixtures of skew-normal
distributions, but very   little is known in the literature for the general case.

The rest of this paper proceeds as follows. In   Section 2, we first provide a   review  of some results on  multivariate elliptical and multivariate skew-normal
distributions.  In   Section 3, we derive    stochastic representations,   characteristic functions  and the first four moments of a skew-elliptical  vector ${\bf Y}$, and
also the first two moments of associated quadratic form ${\bf Y'AY}$ for a symmetric matrix ${\bf A}$. Section 4 presents a   summary of the work carried out here.
Some essential background   details are presented in the Appendix.

\numberwithin{equation}{section}
\section{ Some preliminaries}
An $n$-dimensional random vector $\mathbf{X} = \left(X_{1},...,X_{n}\right)'$ is said to have a multivariate elliptical distribution if its joint density function is given by
 \begin{equation}
 f\left(\mathbf{x}\right)=c_{n}|\mathbf{\Sigma} |^{-\frac{1}{2}}g^{(n)}\left[\mathbf{\left(x-\boldsymbol{\mu}\right)}' \mathbf{\Sigma}^{-1}\mathbf{\left(x-\boldsymbol{\mu}\right)}\right],
 \end{equation}
 where $\boldsymbol{\mu}\in \mathbb{R}^{n}$ is the location vector, $\mathbf{\Sigma}$ is the positive definite  dispersion matrix, $g^{(n)}$ is a non-negative function  satisfying $\int_{0}^{\infty }u^{n/2-1}g^{(n)}\left(u\right)\mathrm{d}u<\infty$, and $c_{n}$ is a normalized constant. In this case, we shall write $\mathbf{X}\sim EC_{n}\left(\boldsymbol{\mu},\mathbf{\Sigma},g^{(n)}\right)$. Here, $g^{(n)}$ is called the density generator.

If $\mathbf{X}\sim EC_{n}\left(\boldsymbol{\mu},\mathbf{\Sigma},g^{(n)}\right)$, then the characteristic function of $\mathbf{X}$ has the form
$$
\varphi\left(\mathbf{t}\right) = \exp\left(i\mathbf{t}'\boldsymbol{\mu}\right)\psi\left(\mathbf{t}'\mathbf{\Sigma}\mathbf{t}\right),~~~\mathbf{t}\in \mathbb{R}^{n},
$$
where $\psi$ is a real-valued function with $\psi(0)=1$. The notation $\mathbf{X}\sim EC_{n}\left(\boldsymbol{\mu},\mathbf{\Sigma},\psi\right)$ is also used
 for the $n$-dimensional elliptical distributions generated from the
function  $\psi$.  The covariance matrix  of $\mathbf{X}$ is proportional  matrix $\mathbf{\Sigma}$, i.e. $Cov\left(\mathbf{X}\right)=-2\psi'\left(0\right)\mathbf{\Sigma}$.

 We start with a conclusion from Fang et al. (1990). The random vector  $\mathbf{X}\sim EC_{n}\left(\boldsymbol{\mu},\mathbf{\Sigma},g^{(n)}\right)$ if and only if
  ${\bf X}\overset{d}{=} \boldsymbol{\mu}+ R  \mathbf{\Sigma}^{\frac12} {\bf U}^{(n)},$
where $R\ge  0$ is a random variable independent of ${\bf U}^{(n)}$, $R$ is one to one with the generator function $g$, $\mathbf{\Sigma}^{\frac12}$  is the  unique positive definite square root matrix of $\mathbf{\Sigma}$,   and ${\bf U}^{(n)}$   is the uniform random
vector on the unit sphere ${\cal S}^{n-1}=\{{\bf u}\in \Bbb{R}^n: ||{\bf u}||^2=1\} $.
 Moreover,
$({\bf X}-\boldsymbol{\mu})'\mathbf{\Sigma}^{-1}({\bf X}-\boldsymbol{\mu})\overset{d}{=}R^2$.  The notation ${\bf X}\overset{d}{=}{\bf Y}$
  means that  ${\bf X}$ and ${\bf Y}$ have the same distribution.

An $n$-dimensional random (column) vector ${\bf X}$  is said to have a spherical (symmetric) distribution if its characteristic function has the form
$\phi_{\bf X}({\bf t})=\phi({\bf t}'{\bf t}), {\bf t}\in \Bbb{R}^n$, for
some function $\phi$,  and we denote it by ${\bf X}\sim S_n(\phi)$.
Cambanis et al. (1981) established the fact that   ${\bf X}\sim S_n(\phi)$ if and only if $ {\bf X}$ has a stochastic representation ${\bf X}\overset{d}{=}R {\bf U}^{(n)},$
 where $R=||{\bf X}||$  is independent of  ${\bf U}^{(n)}$ and the distribution function $F_R$ of $R$ is related
to $\phi$ through
$\phi(u)=\int_0^{\infty}\Omega_n(r^2 u)dF_R(r)$,
where $\Omega_n({\bf t}'{\bf t})=E(e^{i{\bf t}'{\bf U}^{(n)}})$ is the characteristic function of ${\bf U}^{(n)}$.    Moreover,
 ${\bf X}$ possesses a density generator $g$ if and
only if $R$ has a density $f_R$, and the relationship between $g$ and $f_R$ is as follows (see, e.g. Fang et al. (1990)):
  $$f_R(r)=\frac{2\pi^{n/2}}{\Gamma(n/2)}r^{n-1}g(r^2),\; r>0.$$
 When a random vector  ${\bf X}$ is uniformly distributed  inside  the sphere ${\cal S}^{n-1}$ in $\Bbb{R}^n$,
the random variable $R$  has a  Beta$(n, 1)$ distribution. When a random vector   ${\bf X}\sim N_n({\bf 0}, {\bf I}_n)$,  the random variable $R$ has a $\chi$-distribution with $n$ degrees of freedom.
 For more information about some other subclasses of $n$-dimensional  elliptically contoured distributions,  interested readers may refer to    Fang et al. (1990).

Let  ${\bf U}^{(n)}=(U_1, \dots, U_n)'$. Then, for each $1\le k\le n-1$,   $(U_1,\dots, U_k)'$ has a  multivariate
Pearson type II distribution, $PII_k({\bf 0},{\bf I}_k,\frac{n-k-2}{2})$,   with   density function [see, for instance, Stam (1982)]
$$p(u_1,\cdots,u_k)=\frac{\Gamma(\frac{n}{2})}{\pi^{\frac{k}{2}}\Gamma(\frac{n-k}{2})}\left(1-\sum_{i=1}^{k}u_i^2\right)^{\frac{n-k-2}{2}}, \;\; 0<\sum_{i=1}^{k}u_i^2<1.$$
In particular,  $U_1$ has   a univariate Pearson type II distribution
with p.d.f.
$$p(u_1)=\frac{\Gamma(\frac{n}{2})}{\pi^{\frac{1}{2}}\Gamma(\frac{n-1}{2})}\left(1-u_1^2\right)^{\frac{n-3}{2}},  \;\;  -1<u_1<1.$$
The following  important stochastic representation for ${\bf X}\sim S_n(\phi)$, from  Cambanis,  Huang and Simons (1981), will be used repeatedly  in the sequel.
 \begin{lemma}  (Cambanis et al. (1981)). If   ${\bf X}={({{\bf X}_1}', {\bf X_2}')}'\overset{d}{=}R {\bf U}^{(n)}$
 $\sim~S_n(\phi)$, where    ${\bf X}_1$ and $ {\bf X}_2$  have $m$ ($1\le m\le n-1$) and $n-m$
components, respectively, then
 \begin{eqnarray*}
\left(\begin{array}{cc}
 {\bf X}_1&\\
 {\bf X}_2&
\end{array}\right)
=\left(
\begin{array}{cc}
 &R_1{\bf U}^{(m)}\\
 &R_2{\bf U}^{(n-m)}
\end{array}
\right),
 \end{eqnarray*}
 where $R_1 \overset{d}{=} Rd_1$ and $ R_2 \overset{d}{=} Rd_2$. Here, $R_1=||{\bf X}_1||=\sqrt{X_1^2+\cdots+X_m^2}, R_2=\sqrt{X_{m+1}^2+\cdots+X_n^2}, d_1^2=\frac{||{\bf X}_1||^2}{||{\bf X}||^2}\sim Beta(\frac{m}{2},\frac{n-m}{2})$ and $d_2^2=1-d_1^2$.
  Moreover, $R, (d_1,d_2), {\bf U}^{(m)},  {\bf U}^{(n-m)}$ are all independent.
\end{lemma}
If $\Phi_k, k\ge 1$, is the class of all functions $\phi: [0, \infty)\rightarrow \Bbb{R}$ such
that $\phi({\bf t}'{\bf t}), {\bf t}\in {\Bbb R}^k$, is a characteristic function, then  $\phi\in \Phi_k$  if and only if
$$\phi(u)=\int_0^{\infty}\Omega_k(r^2 u)dF_k(r),$$
where $\Omega_k({\bf t}'{\bf t})=E(e^{i{\bf t}'{\bf U}^{(k)}})$ is the characteristic function of ${\bf U}^{(k)}$.  Furthermore, if  $\phi\in \Phi_k$, then for each $1\le m\le k$, we have $\phi\in  \Phi_m$ and
$$\phi(u)=\int_0^{\infty}\Omega_m(r^2 u)dF_m(r).$$
The precise relationship between the various $F_n$'s, for $1\le n\le k$, is as  in the following lemma.
\begin{lemma}  (Cambanis et al. (1981)).  If $1\le m<n\le k$,  $R_m\sim F_m$ and  $R_n\sim F_n$, then  $R_m \overset{d}{=} R_{mn}R_n,$ where
$R_{mn}^2\sim Beta(m/2, (n-m)/2).$
\end{lemma}

\numberwithin{equation}{section}
\section{Skew-elliptical distributions }\label{intro}

The skew-elliptical  (SE) family of distributions is a  flexible class of location-scale models. Extensive discussions on this family of distributions have been provided by Branco and Dey (2001), Azzalini and Capitanio (1999, 2003), Fang (2003, 2005),  Arellano-Valle
and Azzalini (2006), Arellano-Valle and Genton (2005, 2010a, 2010b), and  Genton (2004). This family of distributions  facilitates capturing  the skewness and kurtosis present in the distribution of the data.
Recent developments relating  to multivariate skew-normal and skew-elliptical  distributions have been reviewed by   Azzalini (2022) and Lee and McLachlan (2022). In this section, we present  some stochastic representations, characteristic functions and moments of  skew elliptical distributions.

\subsection{Stochastic representations  }

Stochastic representations are useful for random number generation, and additionally, they may   provide a motivation for the adoption
of skew-ellptical family as a stochastic model for observed data.
A direct approach  for defining a SE distribution  is  through  conditioning mechanism.  We first   describe  the
SE distribution  along    the constructive route of   Branco and Dey (2001), which extended Azzalini and
Dalla Valle's (1996) results to multivariate  SE distributions.

  An $n$-dimensional absolutely continuous random vector ${\bf X}$ is said to have skew-elliptical distribution with location vector $\boldsymbol{\mu}\in \Bbb{R}^n$,
 positive definite dispersion matrix $\bf{\Omega}$
 of dimension $n\times n$, skewness/shape vector $\boldsymbol{\alpha}\in \Bbb{R}^n$ and density generating function $g^{(n+1)}$, in symbols
${\bf X}\sim SE_{n}\left(\boldsymbol{\mu},\mathbf{\Omega},\boldsymbol{\alpha}, g^{(n+1)}\right)$,  if its density function is given
by
\begin{equation}\label{eq12}
f_{\mathbf{X}}(\mathbf{x})=2f_{g^{(n)}}(\mathbf{x};\boldsymbol{\mu},\mathbf{\Omega})F_{g_{q(\mathbf{x})}}\left(\boldsymbol{\alpha}'\boldsymbol{\omega}^{-1}(\mathbf{x}-\boldsymbol{\mu})\right),\; \mathbf{x}\in \Bbb{R}^n,
 \end{equation}
 where $f_{g^{(n)}}(\cdot;)$ is the density  of   $EC_{n}\left(\boldsymbol{\mu},\mathbf{\Omega},g^{(n)}\right)$, and $F_{g_{q(\mathbf{y})}}(\cdot)$ is the cumulative distribution function of  $EC_{1}\left(0,1, g_{q(\mathbf{y})}\right)$. Here, $\boldsymbol {\omega}={\rm diag}(\mathbf{\Omega})^{\frac12}$, and
 \begin{eqnarray*}
& \boldsymbol{\alpha}=\mathbf{\Omega}^{-1}\boldsymbol{\delta}\left(1-\boldsymbol{\delta}^{T}\mathbf{\Omega}^{-1}\boldsymbol{\delta}\right)^{-\frac12}, \;
g^{(n)}(u)=\int_{0}^{\infty}r^{-\frac{1}{2}}g^{(n+1)}(r+u)\mathrm{d}r, \\
&g_{q(\mathbf{y})}(u)=g^{(n+1)}\left[u+q(\mathbf{y})\right]/g^{(n)}\left(q(\mathbf{y})\right), \;
q(\mathbf{y})=(\mathbf{y}-\boldsymbol{\mu})^{T}\mathbf{\Omega}^{-1}(\mathbf{y}-\boldsymbol{\mu}). 	
\end{eqnarray*}
Branco and Dey (2001) adopted the notation  $SE_{n}\left(\boldsymbol{\mu},\mathbf{\Omega},\boldsymbol{\delta}, g^{(n+1)}\right)$ for this distribution.
When  $\boldsymbol{\alpha}={\bf 0}$, the density in (3.1) reduces to the  multivariate  elliptical distribution $EC_{n}\left(\boldsymbol{\mu},\mathbf{\Omega},g^{(n)}\right)$.

 The class of skew elliptical densities in (3.1)  obtained through the conditioning method is equivalent to the one obtained by applying the transformation method.  In
Proposition 9 of  Azzalini and Capitanio (2003)  (see also Fang (2003)) , it has been shown that, as in the case of   skew-normal family, the skew-elliptical vector   ${\bf Y}\sim SE_{n}\left(\boldsymbol{0},\mathbf{{\bar\Omega}},\boldsymbol{\alpha}, g^{(n+1)}\right)$ can be obtained by an additive
construction,   where $\mathbf{\bar\Omega}=  \boldsymbol{\omega}^{-1}\mathbf{\Omega} \boldsymbol{\omega}^{-1}$.  In fact, we
 consider the   $(n+1)$-dimensional random vector
 \begin{eqnarray}
\left(\begin{array}{cc}
 {U}_0&\\
 {\bf U}&
\end{array}
\right)
\sim EC_{1+n}\left(
\left(\begin{array}{cc}
 &0\\
 &{\bf 0}
\end{array}\right),
 \left(
\begin{array}{ccc}
  1 &  \boldsymbol{0}'  \\
 \boldsymbol{0}  & \mathbf{\Psi}
\end{array}
\right),
 g^{(n+1)}
\right),
 \end{eqnarray}
 where $\boldsymbol{0}=(0_{1},0_{2},...,0_{n})'$  and
    $\mathbf{\Psi}$ is a $n\times n$ full rank correlation matrix. Let us define
 \begin{eqnarray}
 {\bf Y}=\boldsymbol{\delta} |U_0|+  \mathbf{\Delta}{\bf U},
  \end{eqnarray}
  where    $\boldsymbol{\delta}=(\delta_1,\cdots,\delta_n)'$ with $\delta_1,\cdots,\delta_n\in (-1,1)$, and
  $$ \mathbf{\Delta}={\rm diag}\left\{\left(1-\delta_{1}^{2}\right)^{\frac{1}{2}},\cdots,\left(1-\delta_{n}^{2}\right)^{\frac{1}{2}}\right\}.$$
   Then, the  density of  ${\bf Y}$ is of the form
   \begin{equation}\label{eq12}
f_{\mathbf{Y}}(\mathbf{y})=2f_{g^{(n)}}(\mathbf{y}; \boldsymbol{0},\mathbf{\bar{\Omega}})F_{g_{q(\mathbf{y})}}\left(\boldsymbol{\alpha}^{T}\mathbf{y}\right),\; \mathbf{y}\in \Bbb{R}^n.
 \end{equation}
  The   relationship   between  $(\mathbf{\Psi}, \boldsymbol{\delta})$ and  $(\mathbf{\bar{\Omega}},  \boldsymbol{\alpha})$ are then
   \begin{eqnarray*}
& \mathbf{\bar{\Omega}}&=\mathbf{\Delta}(\mathbf{\Psi}+ \boldsymbol{\lambda} \boldsymbol{\lambda}')\mathbf{\Delta}
=\mathbf{\Delta}\mathbf{\Psi}\mathbf{\Delta}+\boldsymbol{\delta} \boldsymbol{\delta}', \\
& \boldsymbol{\alpha}&=(1+ \boldsymbol{\lambda}'\mathbf{\Psi}^{-1} \boldsymbol{\lambda})^{-\frac12}\mathbf{\Delta}^{-1}\mathbf{\Psi}^{-1} \boldsymbol{\lambda}, 	
\end{eqnarray*}
where $\boldsymbol{\lambda}=\mathbf{\Delta}^{-1}  \boldsymbol{\delta}.$
To incorporate location and scale parameters, we suppose that   $\boldsymbol{\mu}=(\mu_1,\cdots,\mu_n)'$ and
 $\mathbf{\Omega}=(w_{ij})$ is a full rank $n\times n$ covariance matrix of ${\bf Y}$.
    Then, the density function of the location-scale vector  ${\bf X}:=\boldsymbol{\mu}+ \boldsymbol{\omega}{\bf Y}$, with ${\bf Y}$ defined by (3.2) and   $\boldsymbol {\omega}={\rm diag}(\mathbf{\Omega})^{\frac12}$,
     is given by (3.1). See     Azzalini and Capitanio  (2003) for more details.

We present the above classical results in the form of the following lemma.
  \begin{lemma}
If ${\bf X}\sim SE_{n}\left(\boldsymbol{\mu}, \mathbf{ \Omega},\boldsymbol{\alpha}, g^{(n+1)}\right)$,  then  ${\bf X}$ admits the
stochastic representation
 $${\bf X} \stackrel{d}{=} \boldsymbol{\omega}(\boldsymbol{\delta} |U_0|+  \mathbf{\Delta}{\bf U})+ \boldsymbol{\mu},$$
 where $U_0$ and  ${\bf U}$ have a joint elliptical distribution given by (3.2).
\end{lemma}
In the sequel,   based  on  Lemma 3.1, we formulate several  additional stochastic representations for skew-elliptical random vectors.
\begin{theorem}
Suppose ${\bf X}\sim SE_{n}\left(\boldsymbol{\mu}, \mathbf{ \Omega},\boldsymbol{\alpha}, g^{(n+1)}\right)$. Then,   ${\bf X}$ admits the following equivalent
stochastic representations:\\
$(a) \; \; {\bf X}\stackrel{d}{=}R_0\left( \boldsymbol{\omega}\boldsymbol{\delta} |U_1|+ \boldsymbol{\omega}\mathbf{\Delta} \mathbf{\Psi}^{\frac12}\mathbf{U}_2\right)+\boldsymbol{\mu};\\
(b)\; \; \mathbf{X}\stackrel{d}{=}R_0\left(d_1  \boldsymbol{\omega}\boldsymbol{\delta}|U^{(1)}|+d_2  \boldsymbol{\omega} \mathbf{\Delta}  \mathbf{\Psi}^{\frac12}\mathbf{U}^{(n)}\right)+\boldsymbol{\mu};\\
 (c)\; \; \mathbf{X}\stackrel{d}{=}R_1  \boldsymbol{\omega}\boldsymbol{\delta} |U^{(1)}|+R_n  \boldsymbol{\omega} \mathbf{\Delta} \mathbf{\Psi}^{\frac12}\mathbf{U}^{(n)}+\boldsymbol{\mu},$\\
 where $\mathbf{\Psi}^{\frac12}$ is the unique symmetric square root  of $\mathbf{\Psi}$,  $R_0$ is an absolutely continuous non-negative random variable whose density function  is given by
\begin{equation}
 h_{R_0}(r)=\frac{2\pi^{\frac{n+1}{2}}}{\Gamma\left(\frac{n+1}{2}\right)}r^{n}g^{(n+1)}\left(r^{2}\right), \;r>0,
 \end{equation}
 $\left(U_1,\mathbf{U}_2'\right)'$  is a partition of the
uniform vector ${\bf U}^{(n+1)}$ and independent of $R_0$, $d_1^2\sim Beta(\frac{1}{2},\frac{n-1}{2})$, $d_2^2=1-d_1^2$,  $R_1 \overset{d}{=} d_{1}R_0$,   $R_n \overset{d}{=} \sqrt{1-d^2_1}R_0$, $U_1\stackrel{d}{=}d_1 U^{(1)}$, and ${\bf U}_2 \stackrel{d}{=}d_2 \mathbf{U}^{(n)}$.
  Moreover, $R_0, (d_1,d_2), {U}^{(1)},$  ${\bf U}^{(n)}$ are all independent,  $(R_1, R_n), {U}^{(1)},$  ${\bf U}^{(n)}$ are all independent.
\end{theorem}
{\bf Proof}.  Since ${\bf U}^{(n+1)}=\left(U_1,\mathbf{U}_2'\right)'$
  is the uniform vector on the $(n + 1)$-dimensional unit sphere, then by
(3.2) we have $(U_0, {\bf U})\stackrel{d}{=}(R_0 U_1, R_0 \mathbf{\Psi}^{\frac12} {\bf U}_2)$.
Thus, according to Lemma 3.1, we have that
  $${\bf Y}=\boldsymbol{\omega}^{-1}({\bf X}-\boldsymbol{\mu})\stackrel{d}{=}(\boldsymbol{\delta}, \boldsymbol{\Delta} ) \left(\begin{array}{cc}
 &{U}_0\\
 &{\bf U}
\end{array}
\right)
\stackrel{d}{=}(\boldsymbol{\delta}, \boldsymbol{\Delta} ) \left(\begin{array}{cc}
 &R_0 U_1\\
 &R_0 \mathbf{\Psi}^{\frac12} {\bf U}_2
\end{array}
\right),
$$
which proves part (a). Part (b) follows from (a) and the application of Lemma 2.1, while (c)
is a consequence of (b) and Lemma 2.2. $\hfill\square$

We now give two  examples of multivariate skew-elliptical distributions as special cases.

{\bf Example 3.1} If ${\bf X}\sim SN_{n}\left(\boldsymbol{\mu}, \mathbf{\Omega},\boldsymbol{\alpha}\right)$, then ${\bf X}$ satisfies the stochastic
representations of Theorem 3.1, with  $R_0\sim \chi_{n+1}$, $R_1\sim \chi_1$ and $R_n\sim \chi_n$, where  $\chi_{\nu}^2$
denotes the chi-square distribution with $\nu>0$ degrees of freedom.

 {\bf Example 3.2}  SMSN distribution: An $n$-dimensional random vector ${\bf X}$ is said to have a  scale mixture of skew-normal (SMSN) distributions  if  ${\bf X}$  has the stochastic representation
${\bf X}=\boldsymbol{\mu}+\sqrt{\eta}{\bf Z}$,
 where ${\bf Z}\sim SN_n(\boldsymbol{0},  \mathbf{\Omega},\boldsymbol{\alpha})$
and $\eta$ is a positive random variable with cdf $H(\cdot; \nu)$,  independently of ${\bf Z}$.  Here, $\nu$ is a scalar parameter indexing the distribution of the   mixing  scale factor $\eta$. This distribution is usually denoted by ${\bf X}\sim  SMSN_{n}\left(\boldsymbol{\mu},   \mathbf{\Omega},\boldsymbol{\alpha}, H\right)$. From Theorem 3.1 and Example 3.1,   we can  present several equivalent stochastic representations for ${\bf X}$ in terms of $R_i=\sqrt{\eta}R_{N_i}, i=0,1,n$, where $\eta$ is independent of $R_{N_i}$, and where
$R_{N_0}\sim \chi_{n+1},  R_{N_1}\sim \chi_{1}$  and $R_{N_n}\sim \chi_{n}$.

\numberwithin{equation}{section}
\subsection{Characteristic functions}\label{intro}

The derivation of   characteristic functions or  moment generation functions  for special classes of skewed
distributions has received much attention in the literature;  see, for instance, Azzalini
and Valle (1996),   Nadarajah and Kotz (2003),  Fang (2005),   Kim and Genton (2011),   Vilca,   Balakrishnan and  Zeller (2014)), and  Arellano-Valle and Azzalini (2022).
But, there have been surprisingly few  articles on the derivation of characteristic functions of  general skew-elliptical distributions.
  To the best of our  knowledge, Fang (2003) is the first work in which the  moment generating function of  skew elliptical distribution of
dimension  $k$ has been studied,  giving an expression as an integral on the space of lower dimension, $\min\{4,k+1\} $.
Hence, there is certainly a need to derive   simpler form of  characteristic functions of these skewed distributions.

The following theorem provides one form of the characteristic function of multivariate skew-elliptical distribution.
\begin{theorem}
Suppose   ${\bf X}\sim SE_{n}\left(\boldsymbol{\mu}, \mathbf{ \Omega},\boldsymbol{\alpha}, g^{(n+1)}\right)$ has its p.d.f. as in (3.4).
Then,  its characteristic function  is given by
\begin{eqnarray}
E(e^{i{\bf t}'\mathbf{X}})=2e^{i{\bf t}'\boldsymbol{\mu}} \int_{0}^{\infty} \left(e^{i({\bf t}'\boldsymbol{\delta}_w)u_0}\cdot \phi_{u_0}\left({\bf t}' \boldsymbol{\omega} \mathbf{\Delta} \mathbf{\Psi}\mathbf{\Delta} \boldsymbol{\omega} {\bf t}\right)\right)P(U_0\in du_0),
\end{eqnarray}
where  $\phi_{u_0}(\cdot)$ is a characteristic generator  function, $\boldsymbol{\delta}_w= \boldsymbol{\omega}\boldsymbol{\delta}$ and
$ \boldsymbol{\omega}\mathbf{\Delta}\mathbf{\Psi}\mathbf{\Delta} \boldsymbol{\omega}=\mathbf{\Omega}-\boldsymbol{\delta}_w \boldsymbol{\delta}_w'.$
\end{theorem}
{\bf Proof}. Note that ${\bf X} \stackrel{d}{=}\boldsymbol{\mu}+ \boldsymbol{\omega}{\bf Y}$, where ${\bf Y}$ is as defined in (3.3). For any ${\bf t}\in \Bbb{R}^n$,  applying double expectation
formula, we then get
\begin{eqnarray}
E(e^{i{\bf t}'\mathbf{Y}})&=& E\left(E(e^{i{\bf t}'\mathbf{Y}}|U_0)\right )\nonumber\\
&=&  E\left(E(e^{i({\bf t}'\boldsymbol{\delta})|U_0|})\cdot E(e^{i{\bf t}'(\mathbf{\Delta}\mathbf{U})} |U_0)\right).
\end{eqnarray}
It is well-known that (see Theorem 2.18 of Fang et al. (1990))
$$\mathbf{U}|(U_0=u_0)\sim EC_{n}\left(0,\mathbf{\Psi}, g_{(2.1)}(t)\right),$$
where $g_{(2.1)}(t)=g^{(n)}(t+u_0^2).$
Hence,
$$\mathbf{\Delta}\mathbf{U}|(U_0=u_0)\sim EC_{n}\left(0,\mathbf{\Delta}\mathbf{\Psi}\mathbf{\Delta}, g_{(2.1)}(t)\right).$$
The conditional characteristic function of $\mathbf{\Delta}\mathbf{U}$, given $(U_0=u_0)$,  is of the form
\begin{equation}
   E(e^{i{\bf t}'(\mathbf{\Delta}\mathbf{U})} |U_0=u_0) =\phi_{u_0}\left({\bf t}'\mathbf{\Delta}\mathbf{\Psi}\mathbf{\Delta} {\bf t}\right),
\end{equation}
for some scalar function $\phi_{u_0}(\cdot)$. Substituting (3.9) into (3.8),  we get
\begin{eqnarray*}
E(e^{i{\bf t}'\mathbf{Y}})&=&  \int_{-\infty}^{\infty} \left(E(e^{i({\bf t}'\boldsymbol{\delta})|u_0|})\cdot   E(e^{i{\bf t}'(\mathbf{\Delta}\mathbf{U})} |U_0=u_0) \right)P(U_0\in du_0)\\
&=&  \int_{-\infty}^{\infty} \left(E(e^{i({\bf t}'\boldsymbol{\delta})|u_0|})\cdot \phi_{u_0}\left({\bf t}'\mathbf{\Delta}\mathbf{\Psi}\mathbf{\Delta} {\bf t}\right)\right)P(U_0\in du_0)\\
&=& 2 \int_{0}^{\infty} \left(e^{i({\bf t}'\boldsymbol{\delta})u_0}\cdot \phi_{u_0}\left({\bf t}'\mathbf{\Delta}\mathbf{\Psi}\mathbf{\Delta} {\bf t}\right)\right)P(U_0\in du_0).
\end{eqnarray*}
The result in (3.7) now follows as ${\bf X}=\boldsymbol{\mu}+ \boldsymbol{\omega}{\bf Y}$.  $\hfill\square$

\begin{remark} Letting  $g^{(n)}(u)=(2\pi)^{-n/2}e^{-\frac{u}{2}}$ in Theorem 3.2, we get the  characteristic function  of  ${\bf X}\sim SN_{n}\left(\boldsymbol{\mu}, \mathbf{ \Omega},\boldsymbol{\alpha}\right)$, which has been obtained by Kim and Genton (2011):
 \begin{eqnarray}
 E(e^{i{\bf t}'\mathbf{X}})&=& 2\exp\left(i{\bf t}'\boldsymbol{\mu}-\frac12 {\bf t}'\mathbf{\Omega}{\bf t}\right)\Phi(i\boldsymbol{\delta}'{\bf t})\nonumber\\
 & =&\exp\left( i{\bf t}'\boldsymbol{\mu}-\frac{1}{2}{\bf t}'\mathbf{\Omega} {\bf t}\right)\left\{1+i\tau ( \boldsymbol{{\delta}}' {\bf t}) \right\},\; {\bf t}\in {\Bbb R}^n,\nonumber
   \end{eqnarray}
 where
 $$\boldsymbol{\delta}=\frac{\mathbf{\Omega}\boldsymbol{\alpha}_w}
 {\sqrt{1+\boldsymbol{\alpha}'_w\mathbf{\Omega} \boldsymbol{\alpha}_w}},\;
 \tau(x)=\sqrt{\frac{2}{\pi}}\int_0^x \exp(u^2/2)du,\;\; x>0.$$
\end{remark}

To illustrate the usefulness  of Theorem 3.2 to non-normal distributions,  we consider in the following example  skew-$t$ distribution, which is  an interesting case among scale mixtures of SN distributions.
Kim and Genton (2011) found the characteristic function of multivariate skew-$t$ distribution through  scale mixtures of skew-normal  distributions. We present now
 an alternate expression for this  characteristic function.

\noindent{\bf Example 3.3}. Let ${\bf X}$ follow  skew-$t$ distribution $ St_{n}\left(\boldsymbol{\mu}, \mathbf{ \Omega},\boldsymbol{\alpha}, \nu\right)$,  $(U_0, {\bf U}')'\sim t_{n+1}\left(\boldsymbol{0},\mathbf{\Psi}^*,\nu\right)$, and
\begin{equation*}
\mathbf{\Psi}^* =
\left(
\begin{array}{ccc}
  1 &  \boldsymbol{0}'  \\
 \boldsymbol{0}  & \mathbf{\Psi}
\end{array}
\right),
 \end{equation*}
 where   $\mathbf{\Psi}$ is a $n\times n$ full rank correlation matrix. Then, $U_0\sim t_1(0,1,\nu), {\bf U}\sim t_n({\bf 0}, \mathbf{\Psi},\nu)$ and
 ${\bf U}|U_0\sim t_n\left({\bf 0}, \frac{\nu+U_0^2}{\nu +1}\mathbf{\Psi},\nu+1\right)$ (see Ding (2016)). It then follows that
 $$\Delta{\bf U}|U_0\sim t_n\left({\bf 0}, \frac{\nu+U_0^2}{\nu +1}\Delta\mathbf{\Psi}\Delta,\nu+1\right).$$
 So,
 $$P(U_0\in du_0)=\frac{\Gamma(\frac{\nu+1}{2})}{\sqrt{\pi\nu}\Gamma(\frac{\nu}{2})}\left(1+\frac{u_0^2}{\nu}\right)^{-\frac{\nu+1}{2}}du_0$$
 and
 $$\phi_{L|u_0}\left({\bf t}' \Xi {\bf t}\right)=\frac{||((\nu+u_0^2)  \Xi)^{\frac12}{\bf t}||^{\frac{\nu+1}{2}}}{2^{\frac{\nu+1}{2}-1}\Gamma(\frac{\nu+1}{2})}K_{\frac{\nu+1}{2}}\left(||((\nu+u_0^2)  \Xi)^{\frac12}{\bf t}|| \right),$$
 where $\Xi= \boldsymbol{\omega} \mathbf{\Delta} \mathbf{\Psi}\mathbf{\Delta} \boldsymbol{\omega}$, and $K_{\nu}(\cdot)$ is  the modified Bessel function of the second kind (also called the MacDonald function) (see  Appendix A).   For the characteristic function of
 multivariate $t$-distribution, one may see  Kim and Genton (2011).
 By using formula 4 in Gradshteyn and Ryzhik (2007, p.738), which states that
 \begin{eqnarray}
 &&\int_0^{\infty}(x^2+b^2)^{-\frac12 \nu}\;K_{\nu}(a\sqrt{x^2+b^2})\cos(cx)dx\nonumber\\
 &&=\sqrt{\frac{\pi}{2}}a^{-\nu}b^{\frac12-\nu}(a^2+c^2)^{\frac{\nu}{2}-\frac{1}{4}}
 K_{\nu-\frac12}(b\sqrt{a^2+c^2}),\nonumber
 \end{eqnarray}
 where Re$(a)>0$, Re$(b)>0$ and $c$ is real, we obtain the characteristic function of ${\bf X}$  to be
 \begin{eqnarray}
 &&2\int_{0}^{\infty}\left(e^{i({\bf t}'\boldsymbol{\delta}_w)u_0}\cdot \phi_{L|u_0}\left({\bf t}'  \Xi {\bf t}\right)\right)P(U_0\in du_0)\nonumber\\
 &=& \frac{2\nu ^{\frac{\nu}{2}} ({\bf t}'  \Xi {\bf t})^{\frac{\nu}{2}}}{2^{\frac{\nu-1}{2}} \sqrt{\pi}\Gamma(\frac{\nu}{2})}\int_{0}^{\infty} e^{i({\bf t}'\boldsymbol{\delta}_w)u_0}(\nu+u_0^2)^{-\frac{\nu+1}{4}}K_{\frac{\nu+1}{2}}\left(\sqrt{{\bf t}'  \Xi {\bf t}}
 \sqrt{\nu+u_0^2}\right)du_0\nonumber\\
 &=&\frac{(\nu{\bf t}'\Omega {\bf t})^{\frac{\nu}{4}}}{2^{\frac{\nu}{2}-1} \Gamma(\frac{\nu}{2})} K_{\frac{\nu}{2}}\left(\sqrt{\nu {\bf t}'  \Omega {\bf t}}
 \right)\nonumber\\
 &&+ \frac{2\nu ^{\frac{\nu}{2}} ({\bf t}'  \Xi {\bf t})^{\frac{\nu}{2}}}{2^{\frac{\nu-1}{2}} \sqrt{\pi}\Gamma(\frac{\nu}{2})}i\int_{0}^{\infty} \sin({\bf t}'\boldsymbol{\delta}_w u_0)(\nu+u_0^2)^{-\frac{\nu+1}{4}}K_{\frac{\nu+1}{2}}\left(\sqrt{{\bf t}'  \Xi {\bf t}}
 \sqrt{\nu+u_0^2}\right)du_0,\nonumber
 \end{eqnarray}
 wherein we have used the fact that $ e^{i({\bf t}'\boldsymbol{\delta}_w)u_0}=\cos({\bf t}'\boldsymbol{\delta}_w u_0 )+i\sin( {\bf t}'\boldsymbol{\delta}_w u_0)$. In particular, when $\nu = 1$, we get the characteristic function of  multivariate skew-Cauchy distribution;
when $\nu\to\infty$, we  naturally  deduce the characteristic function of  multivariate skew-normal distribution.

\begin{remark} The characteristic function of   $\mathbf{X}\sim SE_{n}\left( \boldsymbol{\mu},\mathbf{\Omega},\boldsymbol{\alpha}, g^{(n+1)}\right)$
 can be written in the product form
 $$e^{i{\bf t}'\boldsymbol{\mu}}E\left(e^{i({\bf t}'\boldsymbol{\delta}_w)|U_0|}\right) \phi\left({\bf t}' \mathbf{\Omega}{\bf t}-{\bf t}'\boldsymbol{\delta}_w \boldsymbol{\delta}_w'{\bf t}\right)$$
 if and only if   $\mathbf{X}\sim SN_{n}\left( \boldsymbol{\mu},\mathbf{\Omega},\boldsymbol{\alpha}\right)$.    This result can be reduced from (3.7) and a  result, due to Kelker (1970),  that any multivariate elliptical
distribution with mutually independent components must necessarily be multivariate normal.  This provides a negative answer to the  conjecture made by   Genton  and Loperfido (2005)  which states that the characteristic function of generalized skew-elliptical distributions (GSE), or a subclass of GSE distributions, has the  form of (in the case of $\boldsymbol{\mu}={\bf 0}, \mathbf{\Omega}={\bf I}_n$)
 $2 \phi\left({\bf t}'{\bf t}\right)k({\bf t}),$ where the
function $k$ is a function such that $k(-{\bf t}) = 1 - k({\bf t})$  and $0\le k({\bf t})\le 1$.
\end{remark}
 The next result provides another form of  characteristic function of multivariate skew-elliptical distribution.
\begin{theorem}
Suppose   ${\bf X}\sim SE_{n}\left(\boldsymbol{\mu}, \mathbf{ \Omega},\boldsymbol{\alpha}, g^{(n+1)}\right)$ has the stochastic representation in  part (a) of Theorem 3.1.
Then, its characteristic function  is given by
\begin{eqnarray}
E(e^{i{\bf t}'\mathbf{X}})=e^{i{\bf t}'\boldsymbol{\mu}} \int_{0}^{\infty} E(e^{ir_0({\bf t}' \mathbf{M})}) P(R_0\in dr_0), \; {\bf t}\in \Bbb{R}^n,
\end{eqnarray}
where   $\mathbf{M}=\boldsymbol{\delta}_w |U_1|+\mathbf{\Delta} \boldsymbol{\omega}\mathbf{\Psi}^{\frac12}\mathbf{U}_2$.
\end{theorem}
{\bf Proof}. Using  Theorem 3.1,   ${\bf X}$ can be expressed as $\mathbf{X}\stackrel{d}{=}R_0\mathbf{M}+\boldsymbol{\mu}$. It then follows that
\begin{eqnarray}
E(e^{i{\bf t}'\mathbf{X}})&=&e^{i{\bf t}'\boldsymbol{\mu}} E(E(e^{iR_0({\bf t}' \mathbf{M})}|R_0))\nonumber\\
&=&e^{i{\bf t}'\boldsymbol{\mu}} \int_{0}^{\infty} E(e^{ir_0({\bf t}' \mathbf{M})}) P(R_0\in dr_0),\nonumber
\end{eqnarray}
where we have used the independence of $\mathbf{M}$ and $R_0$. This completes the proof.  $\hfill\square$

 Recall that  ${\bf U}^{(n)}$ denotes a random vector that is uniformly distributed on the unit sphere ${\cal S}^{n-1}=\{{\bf u}\in \Bbb{R}^n: ||{\bf u}||^2:={\bf u}'{\bf u}=1\}$ and  $\Omega_n(||{\bf t}||^2), {\bf t}\in \Bbb{R}^n $, denotes the  characteristic function of ${\bf U}^{(n)}$. The following results can be found in  Fang et al. (1990).
\begin{lemma} The  characteristic function  $\Omega_n(||{\bf t}||^2)\; ({\bf t}\in \Bbb{R}^n)$ can be expressed in the following equivalent forms:
 \begin{eqnarray}
\Omega_n(||{\bf t}||^2)&=&\Gamma\left(\frac{n}{2}\right)\left(\frac{2}{||{\bf t}||}\right)^{\frac{n-2}{2}}J_{\frac{n-2}{2}}(||{\bf t}||),\\
 \Omega_n(||{\bf t}||^2)&=&\frac{\Gamma\left(\frac{n}{2}\right)}{\Gamma\left(\frac{n-1}{2}\right)\sqrt{\pi}}\int_{-1}^{1} e^{i||{\bf t}||u}(1-u^2)^{\frac{n-1}{2}-1}du,\\
 \Omega_{n}(\|{\bf t}\|^{2})&=&\frac{1}{B(\frac{n-1}{2},\frac{1}{2})}\int_{0}^{\pi}\exp(i\|{\bf t}\|^{2}\cos\theta)\sin^{n-2}\theta d\theta,\\
  \Omega_{n}(\|{\bf t}\|^{2})&=&\frac{\Gamma(\frac{n}{2})}{\sqrt{\pi}}\sum_{k=0}^{\infty}\frac{(-1)^{k}\|{\bf t}\|^{2k}}{(2k)!}\frac{\Gamma(\frac{2k+1}{2})}{\Gamma(\frac{n+2k}{2})},\\
  \Omega_{n}(\|{\bf t}\|^{2})&=&_{0}F_{1}\left(\frac{n}{2};-\frac{1}{4}\|{\bf t}\|^{2}\right),
  \end{eqnarray}
where $J_{\nu}$ is the Bessel function of the first kind of order $\nu$ (see Appendix A), $B(\cdot,\cdot)$ is the beta function and
$${}_0F_1(\gamma;z)=\sum_{k=0}^{\infty}\frac{1}{(\gamma)_k}\frac{z^k}{k!}$$
is the generalized hypergeometric series of order $(0, 1)$.
\end{lemma}

In the following, we derive several equivalent expressions for the characteristic function  of  skew-elliptical distributions.  In order to simplify the presentation of some formulas, we  introduce the notation
$$\Xi(R_1)=\frac{e^{i R_1({\bf t}'\boldsymbol{\delta}_w)}-1}{R_1}.$$
\begin{theorem}
If  ${\bf X}\sim SE_{n}\left(\boldsymbol{\mu}, \mathbf{ \Omega},\boldsymbol{\alpha}, g^{(n+1)}\right)$,  then
\begin{eqnarray}
E(e^{i{\bf t}'\mathbf{X}})=\frac{e^{i{\bf t}'\boldsymbol{\mu}}}{i{\bf t}'\boldsymbol{\delta}_w}E\left( \Xi(R_1) \Omega_n\left(R_n^2||\mathbf{\Delta}
\mathbf{\Omega}^{\frac{1}{2}} \boldsymbol{\omega}{\bf t}||^2\right)\right),\;\; {\bf t}\in \Bbb{R}^n,
\end{eqnarray}
where  $ R_1=R_0d_1$ and $R_n=R_0\sqrt{1-d_1^2}$.  Here, $d_1^2\sim Beta(\frac{1}{2},\frac{n-1}{2})$ and  $R_0$ is a non-negative random variable with its density as in (3.6). Moreover, $R_0$ and $d_1$ are independent.
\end{theorem}
{\bf Proof}.  For any ${\bf t}\in \Bbb{R}^n$, from Theorem 3.1, we have
\begin{eqnarray*}
E(e^{i{\bf t}'\mathbf{X}})&=&E\left(E(e^{i{\bf t}'\mathbf{X}}|R_0, d_1)\right )\\
&=&e^{i{\bf t}'\boldsymbol{\mu}} E\left(E(e^{iR_0d_1({\bf t}'\boldsymbol{\delta}_w)|U^{(1)}|} |R_0, d_1)\cdot E(e^{iR_0\sqrt{1-d^2_1}({\bf t}'
 \boldsymbol{\omega}\mathbf{\Omega}^{\frac{1}{2}}\mathbf{\Delta})\mathbf{U}^{(n)}} |R_0, d_1)\right).
\end{eqnarray*}
Note that $|U^{(1)}|\sim U[0,1]$, and so
$$E(e^{il|U^{(1)}|})=\frac{e^{il}-1}{il}.$$
It then follows that
\begin{eqnarray}
E(e^{iR_0d_1({\bf t}'\boldsymbol{\delta}_w)|U^{(1)}|} |R_0, d_1) =\frac{e^{i R_0d_1({\bf t}' \boldsymbol{\delta}_w)}-1}{i R_0d_1({\bf t}' \boldsymbol{\delta}_w)}.
\end{eqnarray}
Similarly,
\begin{eqnarray}
E(e^{iR_0\sqrt{1-d^2_1}({\bf t}' \boldsymbol{\omega}\mathbf{\Omega}^{\frac{1}{2}}\mathbf{\Delta})\mathbf{U}^{(n)}} |R_0, d_1)
=\Omega_n\left(R_n^2||\mathbf{\Delta}\mathbf{\Omega}^{\frac{1}{2}}  \boldsymbol{\omega}{\bf t}||^2\right).
\end{eqnarray}
The result in (3.16) now follows from (3.17) and (3.18), thus completing   the proof of the theorem.  $\hfill\square$

Using Lemma 3.2 and Theorem 3.4, we readily obtain the following corollary.
\begin{corollary}
  Suppose $\mathbf{X}\sim SE_{n}\left( \boldsymbol{\mu},\mathbf{\Omega},\boldsymbol{\alpha}, g^{(n+1)}\right)$.
   Then, the characteristic function $\varphi_{\mathbf{X}}({\bf t})$ $ ({\bf t}\in \Bbb{R}^n)$ of  $\mathbf{X}$ can be expressed in the following equivalent forms:
 \begin{eqnarray*}
\varphi_{\mathbf{X}}({\bf t})&=&\frac{\Gamma\left(\frac{n}{2}\right)2^{\frac{n-2}{2}} e^{i{\bf t}'\boldsymbol{\mu}}}{i{\bf t}' \boldsymbol{\delta}_w}E\left(\Xi(R_1) \left(\frac{1}{ R_n ||\mathbf{\Delta}\mathbf{\Omega}^{\frac{1}{2}}  \boldsymbol{\omega}{\bf t}||}\right)^{\frac{n-2}{2}}J_{\frac{n-2}{2}}(R_n ||\mathbf{\Delta}\mathbf{\Omega}^{\frac{1}{2}}  \boldsymbol{\omega}{\bf t}||\right),\\
\varphi_{\mathbf{X}}({\bf t})&=&\frac{\Gamma\left(\frac{n}{2}\right)}{\Gamma\left(\frac{n-1}{2}\right)\sqrt{\pi}}\frac{ e^{i{\bf t}'\boldsymbol{\mu}}}{i{\bf t}' \boldsymbol{\delta}_w}E\left(\Xi(R_1) \int_{-1}^{1} e^{i R_n ||\mathbf{\Delta}\mathbf{\Omega}^{\frac{1}{2}} \boldsymbol{\omega} {\bf t}||u}(1-u^2)^{\frac{n-1}{2}-1}du \right),\\
\varphi_{\mathbf{X}}({\bf t})&=&\frac{1}{B(\frac{n-1}{2},\frac{1}{2})}\frac{e^{i{\bf t}'\boldsymbol{\mu}}}{i{\bf t}'\boldsymbol{\delta}_w}E\left(\Xi(R_1)  \int_{0}^{\pi}\exp(iR_n^2||\mathbf{\Delta}\mathbf{\Omega}^{\frac{1}{2}} \boldsymbol{\omega}{\bf t}||^2 \cos\theta)\sin^{n-2}\theta d\theta\right),\\
\varphi_{\mathbf{X}}({\bf t})& =&\frac{\Gamma(\frac{n}{2})}{\sqrt{\pi}}\frac{e^{i{\bf t}'\boldsymbol{\mu}}}{i{\bf t}'\boldsymbol{\delta}_w}E\left(\Xi(R_1)   \sum_{k=0}^{\infty}\frac{(-1)^{k} R_n^2||\mathbf{\Delta}\mathbf{\Omega}^{\frac{1}{2}}  \boldsymbol{\omega}{\bf t}||^{2k}}{(2k)!}
 \frac{\Gamma(\frac{2k+1}{2})}{\Gamma(\frac{n+2k}{2})}\right),\\
\varphi_{\mathbf{X}}({\bf t}) &=&\frac{e^{i{\bf t}'\boldsymbol{\mu}}}{i{\bf t}' \boldsymbol{\delta}_w}E\left(\Xi(R_1)   {}_0F_1\left(\frac{n}{2}; -\frac{R_n^2||\mathbf{\Delta}\mathbf{\Omega}^{\frac{1}{2}}  \boldsymbol{\omega}{\bf t}||^2}{4}\right)\right).
\end{eqnarray*}
\end{corollary}
\begin{corollary}
   Let  $\mathbf{Y}\sim SU^{(n)}\left(\boldsymbol{\delta}\right)$  be the skew-uniform vector with
stochastic representation
\begin{eqnarray*}\mathbf{Y}\stackrel{d}{=} d_1  \boldsymbol{\delta}|U^{(1)}|+d_2  \mathbf{\Delta} \mathbf{U}^{(n)},
 \end{eqnarray*}
where    $d_1^2\sim Beta(\frac{1}{2},\frac{n-1}{2})$ and $d_2^2=1-d_1^2$.
  Moreover, $(d_1,d_2), {U}^{(1)},$  ${\bf U}^{(n)}$ all are  independent.
Then
\begin{eqnarray}
E(e^{i{\bf t}'\mathbf{Y}})&=&\frac{1}{i({\bf t}'\boldsymbol{\delta})}E\left(\frac{e^{i d_1({\bf t}'\boldsymbol{\delta})}-1}{d_1}\Omega_n((1-d_1^2)||\mathbf{\Delta}{\bf t}||^2)\right),\\
E(e^{i{\bf t}'\mathbf{Y}})&=&\frac{1}{i({\bf t}'\boldsymbol{\delta})}E\left(\frac{e^{i d_1({\bf t}'\boldsymbol{\delta})}-1}{d_1} {}_0F_1\left(\frac{n}{2}; -\frac{(1-d_1^2)||\mathbf{\Delta}{\bf t}||^2}{4}\right)\right).
\end{eqnarray}
In particular, when $\boldsymbol{\delta}=0$, we get
$$ E(e^{i{\bf t}'\mathbf{Y}})={}_0F_1\left(\frac{n}{2}; -\frac{||{\bf t}||^2}{4}\right),\;\; {\bf t}\in \Bbb{R}^n.$$
\end{corollary}
{\bf Proof}. From Theorem 3.1, we know that $d_1^2\sim Beta(\frac{1}{2},\frac{n-1}{2})$, and so the p.d.f. of  $d_1$ is given by
$$p(x)=\frac{2\Gamma(\frac{n}{2})}{\pi^{\frac{1}{2}}\Gamma(\frac{n-1}{2})}\left(1-x^2\right)^{\frac{n-3}{2}},\;\; 0<x<1.$$
By using  formula  7.512(11) in Gradshteyn and Ryzhik (2007), we obtain
$$ E(e^{i{\bf t}'\mathbf{Y}})=E\left({}_0F_1\left(\frac{n}{2}; -\frac{(1-d_1^2)||{\bf t}||^2}{4}\right)\right)={}_0F_1\left(\frac{n}{2}; -\frac{||{\bf t}||^2}{4}\right), \;{\bf t}\in \Bbb{R}^n,$$
which completes  the proof of the corollary.  $\hfill\square$

\subsection{Moments and quadratic forms}
In this subsection, we shall derive the moments of random vectors with multivariate skew-elliptical distribution and of their quadratic
forms. Moments of skew-normal random vector ${\bf Z}\sim SN_n(\boldsymbol{\mu},\mathbf{\Omega}_z, \boldsymbol{\alpha})$ and its quadratic forms have been studied by Genton et al.
(2001) using the moment generating function derived in Azzalini and Dalla Valle (1996), where $\boldsymbol{\mu}$ and  $\mathbf{\Omega}_z$  are the mean  and correlation matrix of ${\bf Z}$.    Fang (2003) established the formula of  mixed  product moments of multivariate skew-elliptical distribution.
Kim and Kim (2017) extended these results to the case of   ${\bf X}=\boldsymbol{\mu}+ \boldsymbol{\omega}{\bf Z}$   to
incorporate the scale parameter $ \boldsymbol{\omega}$, where   $ \boldsymbol{\omega}={\rm diag}\{\sqrt{w_{11}},\cdots,\sqrt{w_{nn}}\}$. Moments of scale mixtures of skew-normal distributions and of their quadratic forms were then derived using the simple stochastic relationship   that exists between skew-normal distribution and scale mixtures of skew-normal distribution (see Kim and Kim (2017)).
 The moments, up to the fourth order,  for ${\bf Z}\sim SE_n(\boldsymbol{\mu},\mathbf{\Omega}_z, \boldsymbol{\alpha}, g^{(n+1)})$ have been derived by Yu and Yin (2023) by making use of
 the results for normalized skew-normal distribution in Genton et al. (2001).  We now extend these results  to the more general  case of skew-elliptical distributions.
For $\mathbf{Y}\sim SE_{n}\left(\boldsymbol{\mu},\mathbf{\Omega},  \boldsymbol{\alpha}, g^{(n+1)}\right)$, define its moments up to the fourth order by
\begin{eqnarray*}
M_1&=&E({\bf Y}), \;\;M_2=E({\bf Y}\otimes{\bf Y}')= E({\bf Y}{\bf Y}'),\\
M_3&=&E({\bf Y}\otimes{\bf Y}'\otimes{\bf Y})=E(vec({\bf Y}{\bf Y}'){\bf Y}'), \\
M_4&=&E({\bf Y}\otimes{\bf Y}'\otimes{\bf Y}\otimes {\bf Y}')=E(vec({\bf Y}{\bf Y}')vec({\bf Y}{\bf Y}')'),
\end{eqnarray*}
 provided the involved expected values exist, where $\otimes$ is the Kronecker product and $vec$ is the operator that stacks the columns of a matrix into a single vector.
  For the properties of the vec-operator, the commutation matrix, the
Kronecker product and related matrix algebra,   interested readers are referred to
Schott (1997) and Kollo and von Rosen (2005).  We shall now provide  two methods for the  derivation of these moments.
 In order to simplify the presentation of some formulas, we  introduce the following notation:
\begin{eqnarray*}  A_1&=&\boldsymbol{\delta}_w\otimes\boldsymbol{\mu}'\otimes\boldsymbol{\mu}+\boldsymbol{\mu}\otimes\boldsymbol{\delta}_w'\otimes\boldsymbol{\mu} +\boldsymbol{\mu}\otimes\boldsymbol{\mu}'\otimes\boldsymbol{\delta}_w,\\
A_2&=&\mathbf{\Omega}\otimes\boldsymbol{\mu}+\boldsymbol{\mu}\otimes\mathbf{\Omega}+vec(\mathbf{\Omega})\otimes\boldsymbol{\mu}',\\
A_3&=&\boldsymbol{\delta}_w\otimes\mathbf{\Omega}+vec(\mathbf{\Omega})\boldsymbol{\delta}_w'
+\left(\mathbf{I}_{p}\otimes\boldsymbol{\delta}_w\right)\mathbf{\Omega}-\left(\mathbf{I}_{p}\otimes\boldsymbol{\delta}_w\right)
\left(\boldsymbol{\delta}_w\otimes\boldsymbol{\delta}_w'\right),
\end{eqnarray*}
\begin{equation*}
\begin{split}
B_{1}=&\boldsymbol{\delta}_w\otimes\boldsymbol{\mu}'\otimes\boldsymbol{\mu}\otimes\boldsymbol{\mu}'
+\boldsymbol{\mu}\otimes\boldsymbol{\delta}_w'\otimes\boldsymbol{\mu}\otimes\boldsymbol{\mu}'
+\boldsymbol{\mu}\otimes\boldsymbol{\mu}'\otimes\boldsymbol{\delta}_w\otimes\boldsymbol{\mu}'\\
&+\boldsymbol{\mu}\otimes\boldsymbol{\mu}'\otimes\boldsymbol{\mu}\otimes\boldsymbol{\delta}_w',\\
B_{2}=&\mathbf{\Omega}\otimes\boldsymbol{\mu}\otimes\boldsymbol{\mu}'+\boldsymbol{\mu}\otimes\mathbf{\Omega}\otimes\boldsymbol{\mu}'
+vec(\mathbf{\Omega})\otimes\boldsymbol{\mu}'\otimes\boldsymbol{\mu}'+\boldsymbol{\mu}'\otimes\mathbf{\Omega}\otimes\boldsymbol{\mu}\\ &+\boldsymbol{\mu}\otimes\boldsymbol{\mu}\otimes\left(vec(\mathbf{\Omega})\right)'+\boldsymbol{\mu}\otimes\boldsymbol{\mu}'\otimes\mathbf{\Omega},\\
B_{3}=&\boldsymbol{\delta}_w\otimes\mathbf{\Omega}\otimes\boldsymbol{\mu}'+vec(\mathbf{\Omega})\otimes\boldsymbol{\delta}_w'\otimes\boldsymbol{\mu}'
+\left(\left(\mathbf{I}_{p}\otimes\boldsymbol{\delta}_w\right)\mathbf{\Omega}\right)\otimes\boldsymbol{\mu}'+\boldsymbol{\delta}_w'\otimes\mathbf{\Omega}\otimes\boldsymbol{\mu}\\ &+\boldsymbol{\delta}_w\otimes\left(vec(\mathbf{\Omega})\right)'\otimes\boldsymbol{\mu}
+\left(\mathbf{\Omega}\left(\mathbf{I}_{p}\otimes\boldsymbol{\delta}_w'\right)\right)\otimes\boldsymbol{\mu}+\boldsymbol{\mu}'\otimes\boldsymbol{\delta}_w\otimes\mathbf{\Omega}\\
&+\boldsymbol{\mu}'\otimes\left(vec(\mathbf{\Omega})\boldsymbol{\delta}_w'\right)+\boldsymbol{\mu}'\otimes\left(\left(\mathbf{I}_{p}\otimes\boldsymbol{\delta}_w\right)\mathbf{\Omega}
\right)+\boldsymbol{\mu}\otimes\boldsymbol{\delta}_w'\otimes\mathbf{\Omega}+\boldsymbol{\mu}\otimes\boldsymbol{\delta}_w\otimes\left(vec(\mathbf{\Omega})\right)'\\
&+\boldsymbol{\mu}\otimes\left(\mathbf{\Omega}\left(\mathbf{I}_{p}\otimes\boldsymbol{\delta}_w'\right)\right)-\boldsymbol{\delta}_w\otimes\boldsymbol{\delta}_w'
\otimes\boldsymbol{\delta}_w\otimes\boldsymbol{\mu}'-\boldsymbol{\delta}_w'\otimes\boldsymbol{\delta}_w\otimes\boldsymbol{\delta}_w^{T}\otimes\boldsymbol{\mu}\\
&-\boldsymbol{\mu}'\otimes\boldsymbol{\delta}_w\otimes\boldsymbol{\delta}_w'\otimes\boldsymbol{\delta}_w
-\boldsymbol{\mu}\otimes\boldsymbol{\delta}_w'\otimes\boldsymbol{\delta}_w\otimes\boldsymbol{\delta}_w',\\
B_{4}=&\left(\mathbf{I}_{n^{2}}+\mathbf{K}_{n,n}\right)(\mathbf{\Omega}\otimes\mathbf{\Omega})+vec(\mathbf{\Omega})\left(vec(\mathbf{\Omega})\right)'.
\end{split}
\end{equation*}
Here,  $\boldsymbol{\delta}_w={\bf  w}\boldsymbol{\delta}$ and $\mathbf{K}_{n,n}$ is the commutation matrix associated with a $n\times n$ matrix.

 First, we present the following lemma  given  by Kim and Kim (2017).
 \begin{lemma} Suppose $\mathbf{X}\sim SN_{n}\left(\boldsymbol{0},\mathbf{\Omega},  \boldsymbol{\alpha}\right)$.  Then,
\begin {eqnarray*}
(1)\;\; && E({\bf X})=\sqrt{\frac{2}{\pi}}\boldsymbol{\delta}_w,\;\; (2)\;\; E({\bf X}\otimes{\bf X}')=\mathbf{\Omega},\\
(3) \;\; && E({\bf X}\otimes{\bf X}'\otimes{\bf X})=\sqrt{\frac{2}{\pi}}A_3, \;\; (4)\;\;  E({\bf X}\otimes{\bf X}'\otimes{\bf X}\otimes {\bf X}')=B_4.
\end{eqnarray*}
\end{lemma}
\begin{theorem} Suppose $\mathbf{Y}\sim SE_{n}\left(\boldsymbol{\mu},\mathbf{\Omega},  \boldsymbol{\alpha},   g^{(n+1)}\right)$    has the stochastic representation $\mathbf{Y}\stackrel{d}{=}R_y\mathbf{M}+\boldsymbol{\mu}$, with  $\mathbf{M}=\boldsymbol{\delta}_w |U_1|+ {\bf w}\mathbf{\Delta}  \mathbf{\Psi}^{\frac12}\mathbf{U}_2$ and $E(R_y^4)<\infty$.   Then,
\begin{eqnarray*}
M_1 &=& \boldsymbol{\mu}+\sqrt{\frac{2}{\pi}} \frac{E(R_y)}{E(R_{x})}\boldsymbol{\delta}_w,\\
 M_2 &=&\boldsymbol{\mu}\boldsymbol{\mu}'
+\sqrt{\frac{2}{\pi}}\frac{E(R_y)}{E(R_{x})}\left(\boldsymbol{\mu}\boldsymbol{\delta}_w'
+\boldsymbol{\delta}_w \boldsymbol{\mu}'\right)+\frac{E\left(R_y^{2}\right)}{E\left(R_{x}^{2}\right)}\mathbf{\Omega},\\
M_{3}&=&\boldsymbol{\mu}\otimes\boldsymbol{\mu}'\otimes\boldsymbol{\mu}+\sqrt{\frac{2}{\pi}}\frac{E(R_y)}{E(R_{x})} A_1\\
&&+\frac{E\left(R_y^{2}\right)}{E\left(R_{x}^{2}\right)}A_2+\sqrt{\frac{2}{\pi}}\frac{E(R_y^{3})}{E(R_{x}^{3})}A_3, \\
 M_{4}&=&\boldsymbol{\mu}\otimes\boldsymbol{\mu}'\otimes\boldsymbol{\mu}\otimes\boldsymbol{\mu}'+\sqrt{\frac{2}{\pi}}\frac{E(R_y)}{E(R_{x})}B_{1}\\
 &&+\frac{E(R_y^{2})}{E(R_{x}^{2})}B_{2}+\sqrt{\frac{2}{\pi}}\frac{E(R_y^{3})}{E(R_{x}^{3})}B_{3}
+\frac{E(R_y^{4})}{E(R_{x}^{4})}B_{4},
\end{eqnarray*}
where   $\mathbf{X}\sim SN_{n}\left(\boldsymbol{0},\mathbf{\Omega},  \boldsymbol{\alpha}_w\right)$ and  $\mathbf{X}\stackrel{d}{=}R_{x}\mathbf{M}$.
\end{theorem}
{\bf Proof}. We present two methods for establishing the above expressions.
{\bf Method 1}:
By using $\mathbf{X}\stackrel{d}{=}R_{x}\mathbf{M}$  and the independence of $R_x$ and $\mathbf{M}$,    we obtain
\begin{eqnarray*}
&&E(\mathbf{M})=\frac{E(\mathbf{X})}{E(R_{x})},\;\;E\left(\mathbf{M}\mathbf{M}^{'}\right)=\frac{E\left(\mathbf{X}{\mathbf{X}}^{'}\right)}{E\left(R_{x}^{2}\right)},\\
&&E\left(\mathbf{M}\otimes\mathbf{M}^{'}\otimes\mathbf{M}\right)=\frac{E\left(\mathbf{X}\otimes\mathbf{X}^{'}\otimes\mathbf{X}\right)}{E\left(R_{x}^{3}\right)},\\
&&E\left(\mathbf{M}\otimes\mathbf{M}^{'}\otimes\mathbf{M}\otimes\mathbf{M}^{'}\right)=\frac{E\left(\mathbf{X}\otimes\mathbf{X}^{'}\otimes\mathbf{X}\otimes\mathbf{X}^{'}\right)}
{E\left(R_{x}^{4}\right)},
\end{eqnarray*}
  which,  together with  $\mathbf{Y}\stackrel{d}{=}R_y\mathbf{M}+\boldsymbol{\mu}$ and Lemma 3.3, yield the desired results.\\
{\bf  Method 2}:  $M_i$'s    can be found by differentiating the characteristic function of the form in (3.10) by using matrix
derivatives (see  Appendix B). For example, the  characteristic function of $\mathbf{Y}$ is
\begin{eqnarray*}
E(e^{i{\bf t}'\mathbf{Y}})=e^{i{\bf t}'\boldsymbol{\mu}} \int_{0}^{\infty} E(e^{ir({\bf t}' \mathbf{M})}) P(R_y\in dr), \; {\bf t}\in \Bbb{R}^n,
\end{eqnarray*}
where   $\mathbf{M}=\boldsymbol{\delta}_w |U_1|+ \boldsymbol{\omega} \mathbf{\Delta}  \mathbf{\Psi}^{\frac12}\mathbf{U}_2$. Thus,
 \begin{eqnarray*}
M_1 &=&\frac{1}{i}\frac{\partial}{\partial {\bf t}}E(e^{i{\bf t}'\mathbf{Y}}){\big |_{\bf t=0}}
 =\boldsymbol{\mu}+ E(\mathbf{M})\int_0^{\infty}rP(R_y\in d r) \\
&=& \boldsymbol{\mu}+ E(\mathbf{M}) E(R_y)
 =\boldsymbol{\mu}+\sqrt{\frac{2}{\pi}} \frac{E(R_y)}{E(R_{x})}\boldsymbol{\delta}_w.
\end{eqnarray*}   $\hfill\square$
\begin{remark}
Using (3.6), the $k$-th moment of $R_y$ (if it exists) can be written as
\begin{equation*}
E(R_y^k)=\frac{2\pi^{\frac{n+1}{2}}}{\Gamma\left(\frac{n+1}{2}\right)}\int_0^{\infty}r^{n+k}g^{(n+1)}\left(r^{2}\right)dr.
 \end{equation*}
If, in particular,  $g^{(n+1)}(u)=(2\pi)^{-\frac{n+1}{2}}\exp(-\frac{u}{2})$, we obtain the  $k$-th moment of $R_x$  since  $R_x\sim \chi_{n+1}$.
\end{remark}
  Alternatively, we can find $M_i$'s    by differentiating the characteristic function of the form in (3.7) by using matrix
derivatives (see  Appendix B) as follows:
\begin{eqnarray}
M_1 &=&\frac{1}{i}\frac{\partial}{\partial {\bf t}}E(e^{i{\bf t}'\mathbf{Y}}){\big |_{\bf t=0}} =\boldsymbol{\mu}+ E(|U_0|)\boldsymbol{\delta}_w,\\
M_2&=&\frac{1}{i^2}\frac{\partial^2}{\partial {\bf t}\partial {\bf t'}}E(e^{i{\bf t}'\mathbf{Y}}){\big |_{\bf t=0}}\nonumber\\
  &=&\boldsymbol{\mu}\boldsymbol{\mu}'
+ E(|U_0|)\left(\boldsymbol{\mu}\boldsymbol{\delta}_w'
+\boldsymbol{\delta}_w\boldsymbol{\mu}'\right)-2E(\phi_{L|U_0}'(0))\mathbf{\Omega},\\
M_3&=&\frac{1}{i^3}\frac{\partial^3}{\partial {\bf t}\partial {\bf t'}\partial {\bf t}}E(e^{i{\bf t}'\mathbf{Y}}){\big |_{\bf t=0}}\nonumber\\
  &=&\boldsymbol{\mu}\otimes\boldsymbol{\mu}'\otimes\boldsymbol{\mu}+ E(|U_0|)A_1 \nonumber\\
  &&-2E(\phi_{L|U_0}'(0))A_2-2E(\phi_{L|U_0}'(0)|U_0|)A_3 \nonumber \\
&&+[E(|U_0|^3)+4E(\phi_{L|U_0}'(0)|U_0|)](\boldsymbol{\delta}_w\otimes\boldsymbol{\delta}_w'\otimes\boldsymbol{\delta}_w), \\
M_{4}&=&\frac{1}{i^4}\frac{\partial^4}{\partial {\bf t}\partial {\bf t'}\partial {\bf t} \partial {\bf t'}}E(e^{i{\bf t}'\mathbf{Y}}){\big |_{\bf t=0}}\nonumber\\
&=&\boldsymbol{\mu}\otimes\boldsymbol{\mu}'\otimes\boldsymbol{\mu}\otimes\boldsymbol{\mu}'+E(|U_0|)B_{1}-2  E(\phi_{L|U_0}'(0))B_{2}\nonumber\\
&&-2E(\phi_{L|U_0}'(0)|U_0|)|B_{3}+4 E(\phi_{L|U_0}''(0))B_{4},
\end{eqnarray}
where $\phi'(0)$ denotes the first derivative of  $\phi$ at zero.  Comparing (3.21)-(3.24)  with the results of Theorem 3.5,  we find
$$ E(|U_0|)=\frac{E(R_y)}{E(R_{x})}\sqrt{\frac{2}{\pi}},\;\; E(\phi_{L|U_0}'(0))=-\frac{1}{2}\frac{E\left(R_y^{2}\right)}{E\left(R_{x}^{2}\right)},$$
$$E(\phi_{L|U_0}'(0)|U_0|)= -\sqrt{\frac{1}{2\pi}}\frac{E(R_y^{3})}{E(R_{x}^{3})},\; E(\phi_{L|U_0}''(0)) =\frac{1}{4}\frac{E(R_y^{4})}{E(R_{x}^{4})},$$
$$E(|U_0|^3)+4E(\phi_{L|U_0}'(0)|U_0|)=0.$$
\begin{remark} The expressions of the moments derived in Theorem 3.5 can be effectively used to study multivariate measures of skewness of any member of the skew-elliptical family of distributions along the lines of Balakrishnan and Scarpa (2012); see also  Balakrishnan,   Capitanio and  Scarpa (2014).
\end{remark}
Let ${\bf Y}\sim SN_{n}\left(\boldsymbol{\mu},\mathbf{\Omega},  \boldsymbol{\alpha}\right)$, and let $A, B$
be two symmetric $n\times n$ matrices.  Genton et al. (2001)  derived the  moments of its associated quadratic form, such as $E({\bf Y}'{\bf A}{\bf Y})$, $Var({\bf Y}'{\bf A}{\bf Y})$
and $Cov({\bf Y}'{\bf A}{\bf Y},{\bf Y}'{\bf B}{\bf Y})$. With the first four moments of the random vector ${\bf Y}$ given in Theorem 3.5 and  the  relations  $E({\bf Y}'{\bf A}{\bf Y})={\rm tr}(A E({\bf Y}{\bf Y}'))$ and
 $E(({\bf Y}'{\bf A}{\bf Y})({\bf Y}'{\bf B}{\bf Y}))={\rm tr}((A\otimes B) E({\bf Y}\otimes{\bf Y}'\otimes{\bf Y}\otimes {\bf Y}'))$,
 we can readily compute the first two moments of its quadratic form as presented in the following corollary.
\begin{corollary} Suppose $\mathbf{Y}\sim SE_{n}\left(\boldsymbol{\mu},\mathbf{\Omega},  \boldsymbol{\alpha},   g^{(n+1)}\right)$ and $A, B$
are two symmetric $n\times n$ matrices. Then,
\begin{eqnarray*}
&&E({\bf Y}'{\bf A}{\bf Y})={\rm tr}(AM_2)\\
&&\;\;\;\;\;\;\;\;\;\;\;\;\;\;\;\;= \boldsymbol{\mu}'{\bf A}\boldsymbol{\mu}+ 2\sqrt{\frac{2}{\pi}}\frac{ER_y}{ER_x}\boldsymbol{\mu}'{\bf A}\boldsymbol{\delta}_w+\frac{ER_y^2}{ER_x^2}{\rm tr}({\bf A}\mathbf{\Omega}),\\
&&E({\bf Y}'{\bf A}{\bf Y})^2= {\rm tr}((A\otimes A)M_4),\\
&&Var({\bf Y}'{\bf A}{\bf Y})={\rm tr}((A\otimes A)M_4)-({\rm tr}(AM_2))^2,\\
&&Cov({\bf Y}'{\bf A}{\bf Y},{\bf Y}'{\bf B}{\bf Y})={\rm tr}((A\otimes B)M_4)-{\rm tr}(AM_2){\rm tr}(BM_2).
\end{eqnarray*}
\end{corollary}

  \section{Concluding remarks  }
This work establishes  some  properties of multivariate skew-elliptical distributions.  Specifically, the work
  establishes  stochastic representations, expressions of characteristic functions and the calculation of some multivariate moments and moments of quadratic forms.
It is important to mention  that the extension to more general skew models such as multivariate skew-elliptical distributions with additional parameters (cf. Fang (2003, 2005, 2006)), or the  generalized skew-elliptical distributions (cf. Azzalini and Capitanio (2003), and  Genton and Loperfido (2005))  is  a  more complicated and interesting problem.
 In addition, the results presented here  are  purely theoretical, but what is    important will be to make use of the  results   established here  in   statistical inference and model fitting applications.  \\

\noindent{\bf CRediT authorship contribution statement}\\
{\bf C. Yin}: Conceptualization, Methodology,   Writing-original draft. {\bf N. Balakrishnan}: Supervision, Writing - review \& editing, Conceptualization.\\

\noindent{\bf Acknowledgements.}
The authors express their sincere thanks to the Editor and the anonymous reviewers for all their useful comments and suggestions on an earlier version of this manuscript which led to this improved version.
This research   was supported by the National Natural Science Foundation of China (No. 12071251).\\

\noindent {\bf Appendix A.} \; {\bf Special functions}  (see Gradshteyn and Ryzhik (2007)).

Bessel function of the first kind, $J_{\nu}(x)$,  is defined as
$$J_{\nu}(x)=\sum_{k=0}^{\infty}(-1)^k\frac{1}{k!\Gamma(\nu+k+1)}\left(\frac{x}{2}\right)^{\nu+2k},\; x>0,$$
or equivalently
$$J_{\nu}(x)=\frac{x^{\nu}}{2^{\nu}\Gamma(\nu+1)} {}_0F_1\left(\nu+1;-\frac{x^2}{4}\right),\; x>0.$$
Modified Bessel function of the second kind (also called the MacDonald function), $K_{\nu}(x)$,  is defined as
$$K_{\nu}(x)=\left(\frac{2}{x}\right)^{\nu}\frac{\Gamma(\nu+\frac12)}{\sqrt{\pi}}\int_0^{\infty}(1+u^2)^{-(\nu+\frac12)}\cos (xu)du, x>0,\nu>-\frac12.$$
{\bf Appendix B.}\;
{\bf Some  matrix algebra and  matrix   differentiation   } (see  Magnus and Neudecker (1979) and   Kollo and  von Rosen (2005)).

Let $A=(a_{ij})$ be an $m\times n$ matrix and $B$ be an $s\times t$ matrix. Then, the Kronecker product $A\otimes B$ is defined as the $ms\times nt$ matrix
$A\otimes B=(a_{ij}B)$. With  $A_{.j}$ being  the  $j$th column of $A$, we define the $ mn$ column vector vec($A$) as vec($A)=(A_{.1}',\cdots, A_{.n}')'$.
Let $A$ and $B$ be matrices of any order. Then,
$$vec(ABC)=(C'\otimes A)vec(B),\;\; (vec(A'))'vec(B)={\rm tr}(AB).$$
Let $x$, $y$  and $z$ be column vectors of any order. Then,
$$x  \otimes y=vec(yx'), x  \otimes y'=xy'=y' \otimes x,$$
$$x\otimes y'\otimes z=vec(xz')\otimes y',\;\;  vec(x') = vec(x) = x.$$
The commutation matrix $K_{mn}$ is the matrix multiplier that transforms $vec(A)$ to $vec(A')$: For any $m \times n$  matrix $A$, $K_{mn} vec(A)=vec(A')$. Further,
$$\frac{\partial{t'At}}{\partial t}=At+A't;$$
 if $A$ is symmetric, then
$$\frac{\partial{t'At}}{\partial t}=2At,
\frac{\partial{At}}{\partial t}=A', \frac{\partial{At}}{\partial t'}=A, \frac{\partial{t'A}}{\partial t}=A, \frac{\partial{t't}}{\partial t}=2t, \frac{\partial{a't}}{\partial t}=a. $$

\bibliographystyle{model1-num-names}

\end{document}